# Objectivity and Rigor in Classical Italian Algebraic Geometry

*Silvia De Toffoli and Claudio Fontanari*

**ABSTRACT**
The classification of algebraic surfaces by the Italian School of algebraic geometry is universally recognized as a breakthrough in 20th-century mathematics.  The methods by which it was achieved do not, however, meet the modern standard of rigor and therefore appear dubious from a contemporary viewpoint.  In this article, we offer a glimpse into the mathematical practice of the three leading exponents of the Italian School of algebraic geometry: Castelnuovo, Enriques, and Severi.  We then bring into focus their distinctive conception of rigor and intuition.  Unlike what is often assumed today, from their perspective, rigor is neither opposed to intuition nor understood as a unitary phenomenon – Enriques distinguishes between *small-scale rigor* and *large-scale rigor* and Severi between *formal rigor* and *substantial rigor*.  Finally, we turn to the notion of mathematical objectivity.  We draw from our case study in order to advance a multi-dimensional analysis of objectivity.  Specifically, we suggest that various types of rigor may be associated with different conceptions of objectivity: namely objectivity as faithfulness to facts and objectivity as intersubjectivity.

**RÉSUMÉ**
**Objectivité et rigueur dans la géométrie algébrique italienne classique**
La classification des surfaces algébriques par l'école italienne de géométrie algébrique est universellement reconnue comme une avancée majeure en mathématiques au XXe siècle. Les méthodes utilisées ne répondent cependant pas aux normes modernes de rigueur et apparaissent ainsi suspectes d'un point de vue contemporain. Dans cet article, nous commençons par donner un aperçu de la pratique mathématique des trois principaux représentants de l'école italienne de géométrie algébrique : Castelnuovo, Enriques et Severi. Nous nous intéresserons à leur conception particulière de la rigueur et de l'intuition. Contrairement à une idée communément admise aujourd'hui, de leur point de vue, la rigueur n'est ni opposée à l'intuition ni considérée comme un phénomène uniforme—Enriques distingue entre *rigueur à petite échelle* et *rigueur à grande échelle* et Severi entre *rigueur formelle* et *rigueur substantielle*. A partir de cette étude de cas, nous proposons une analyse multidimensionnelle de l'objectivité mathématique. Plus précisément, nous suggérons que différents types de rigueur peuvent être associés à différentes conceptions de l'objectivité, à savoir l'objectivité comme fidélité aux faits et l'objectivité comme intersubjectivité.

**MOTS CLÉS**
Rigueur, Objectivité, Intuition Mathématique, Géométrie Algébrique, Pratique Mathématique

# §1. SEEING IN MATHEMATICS

The classification of algebraic surfaces was a historic breakthrough in mathematics. It was achieved in the span of few decades by the Italian School of algebraic geometry, in which the joint work of Guido Castelnuovo (1865-1952) and Federigo Enriques (1871-1946) played a prominent role. Algebraic surfaces are an especially thorny subject, particularly if compared to algebraic curves:

> A pioneer and one of the founding fathers of algebraic geometry, the German mathematician Max Noether, after seeing the theory of algebraic curves with its elegance, simplicity, and also depth of results, and comparing it with the collection of the existing examples of algebraic surfaces at that time, for which nothing comparable could be found, used to say that algebraic curves were created by God and algebraic surfaces by the Devil.[1]

The study of algebraic surfaces led to the unexpected discovery of "an order of hidden harmonies where a marvelous beauty shines forth" [un ordine di armonie più riposte ove rifulge una meravigliosa bellezza] among the devilish surfaces.[2] The classification was accomplished by adopting a genuinely experimental approach, vividly described by Castelnuovo:

> It is perhaps worth mentioning what was the method we used to find our way in the darkness in which we found ourselves. We built, in an abstract sense, *bien entendu*, a large number of models of surfaces in our space or in higher spaces, and **we distributed these models in two display cases, so to speak**. One contained the regular surfaces for which everything proceeded as it would in the best of all possible worlds; it was possible to carry over to them the most salient properties of plane curves by analogy. But when we tried to verify these properties on the surfaces of the other display case, the irregular ones, trouble began, and exceptions of all kinds appeared. In the end, **the assiduous study of our models had led us to divine some properties that had to exist,** with appropriate modifications, for the surfaces of both display cases; we then put these properties to the test by building new models. If they resisted the test, we sought, in the last phase, a logical justification. With this procedure, which resembles those of experimental sciences, we succeeded in establishing some distinctive traits among the families of surfaces.[3]

---
[1] E. Bombieri and W. Gubler, *Heights in Diophantine Geometry*, Cambridge, Cambridge University Press, 2006, p. xi.
[2] F. Enriques, *Le Superfici Algebriche*, Bologna, Zanichelli, 1949, p. 464.

For the members of the Italian School, exploring algebraic surfaces was a quasi-perceptual affair requiring a form of intellectual (in)sight. It is worth noting that in the decades contemporary with the Italian School's inquiry, certain types of material – mainly plaster – models of algebraic surfaces were widely used for teaching purposes.[4] But even though Italian geometers attributed great importance to intuition and visualization, they "did not use physical models in their research work but preferred to employ the *Gedankenexperiment* [thought experiment]."[5] An even more striking presentation of this attitude – which was quite widespread at the time but which is rare today, if it exists at all – is found in the recollections of one of the former students and then close collaborator of Enriques, the mathematician Fabio Conforto (1909-1954):

> [Enriques] conceived the algebraic world as existing in itself, independently and outside of us [...] When trying to understand this world, we should not strive for an ideal of logical perfection; and least of all should we proceed axiomatically, starting from postulates [...] **The algebraic world exists by itself** [...] **Understanding it is therefore not really a question of correct deduction, but above all a question of "seeing."** Such a conception deeply satisfied Enriques's powerfully intuitive spirit, who often went so far – and in the close company of his students he was pleased with this apparently paradoxical aspect of his thinking – as to not feel the need for a logical demonstration of some property because he "saw." This made him certain of the truth of the proposition in question and satisfied him fully: certainly it prevented him from proceeding further. Having once declared to him that I did not see the truth of a statement that he believed to be evident, but which we had tried in vain to prove logically, he stopped short (we were in the middle of one of our usual walks) and, instead of attempting a last demonstration, he spun his stick, pointing it toward a little dog on a window sill, saying to me: **"Can't you see? For me, it's as if you were telling me you don't see that little dog!"** And yet, that property, which we found a way to include in the volume on rational surfaces, perhaps still awaits a satisfactory demonstration today.[6]

---

[3] G. Castelnuovo, "La Geometria Algebrica e La Scuola Italiana," in *Atti Del Congresso Internazionale Dei Matematici*, Tomo I:191–201, Bologna, 1928, p. 194. Our translation from Italian.
[4] A. Sattelmacher, *Anschauen, Anfassen, Auffassen. Eine Wissensgeschichte Mathematischer Modelle*, Wiesbaden, Springer, 2021.
[5] L. Giacardi, "Models in Mathematical Teaching in Italy (1850-1950)," *Mathematics and Art III,* ESMA 2015: 11–38, 2015, p. 14.

For Enriques, it is not logic that gives us confidence in the truth of some mathematical propositions, but rather a quasi-perceptual ability to visualize. A blind person could be justified in believing the existence of a little dog in front of her by means of a point-by-point description of the visual field to which she does not have direct access. This would involve the description of a two-dimensional image constructed through the specification of a color for each pixel. However, a person without perceptual impairments has access to a better, more direct source of justification: vision. *Mutatis mutandis*, a person lacking the ability to visualize in mathematics could form a justified belief in the truth of a mathematical proposition through logical analysis. The visual mathematician, however, just sees. And her seeing gives her a better, more direct justification.

An analogy may help us differentiate between various proving practices and to elaborate on what "seeing" means in Enriques's sense. According to Riccardo Brasca,[7] a mathematician who writes a proof is like a painter painting a picture.[8] Each of the painter's individual brushstrokes is meaningful and contributes to the whole picture; likewise, each step of the proof is a partial result that contributes to the large-scale architecture of the overall proof. But the story is different for formal proofs, which are mechanically checkable arguments articulated in a specific formal system. Specifically, Brasca has in mind formal proofs that can be designed with the aid of *interactive proof assistants*, such as *Coq* or *Lean*.[9] If mathematicians are like *painters*, proof assistants are like *printers* that create images by specifying a color one pixel at the time. Unlike a painter's brush, a printer's pixel does not have a specific meaning. It is insignificant. The proof assistant, however, needs to successfully carry out the pixel-by-pixel description because the final product must be readable by a machine (or by a mathematician that cannot *see*).[10]

---

[6] F. Conforto, "Federigo Enriques," *Rend. Mat*, 226–52, 1947. Our translation from Italian.
[7] Radio interview, R. Fulci, "La prova del software", Radio3 scienza, 25 June 2021. Available on-line at: https://www.raiplaysound.it/audio/2021/06/La-prova-del-software-c8c16121-e94c-49ca-bbdc-951e30162f5e.html
.

Let us consider different types of proofs from the point of view of their reliability. If *proof* is a success term, as it is generally taken to be, then all proofs are equal with respect to reliability. They all are 100% reliable in tracking mathematical facts. Yet what we have access to are not proofs themselves, but what one of the authors has labeled *simil-proofs*.[11] These are arguments that look like proofs to the relevant agents but that may contain fatal errors and thus fail to be genuine proofs.[12] Whereas comparing the reliability of proofs does not make sense, it is nonetheless useful to compare the reliability of simil-proofs. The more reliable a simil-proof, the more likely it is that it is a genuine proof – in other words, that it is correct. How is reliability achieved? According to a

---

[8] Much earlier, in a letter to Klein, knot theorist Wilhelm Wirtinger also compared mathematicians to painters in the context of a critique of the progressive centrality of abstraction in mathematics. See: M. Epple, "Branch Points of Algebraic Functions and the Beginnings of Modern Knot Theory," *Historia Mathematica,* 22: 371–401, 1995. Thanks to Michael Friedman for this reference.

[9] J. Avigad, "The Mechanization of Mathematics," *Notices of the American Mathematical Society,* 65 (6): 681–90, 2018.

[10] Note that the perceptual ability to see visually is not required to *see* in mathematics. There are multiple examples of blind mathematicians possessing an outstanding ability to see in mathematics. See, for example, A. Jackson, "The World of Blind Mathematicians," *Notices of the American Mathematical Society,* 49 (10): 1246–51, 2002.

[11] S. De Toffoli, "Groundwork for a Fallibilist Account of Mathematics," *The Philosophical Quarterly,* 71 (4): 823–44, 2021, p. 835.

[12] Here is the definition of simil-proofs:

> An argument is a Simil-Proof (SP) when it is shareable, and some agents who have judged all its parts to be correct as a result of checking accept it as a proof. Moreover, the argument broadly satisfies the standards of acceptability of the mathematical community to which it is addressed.

S. De Toffoli, "Groundwork for a Fallibilist Account of Mathematics," art.cit., p. 835.

widely held view in the philosophy of mathematics, it is rigor that underwrites the reliability of simil-proofs.[13] It follows that questioning the reliability of simil-proofs amounts to investigating them in terms of rigor. The question then presents itself: what is rigor, exactly?

According to the received view in the philosophy of mathematics, rigor is associated with the possibility of formalization. A proof is rigorous if it can be converted into a formal proof. Of course, this characterization is general and needs to be unpacked. The idea behind it, however, is simple. With any given proof, it should be possible, with enough time, energy, background knowledge, and logical skill, to produce a formal proof. Note that with a given proof, many distinct formal proofs can be produced since the formalization process may be carried out in a variety of ways. Although mainstream mathematicians generally do not bother formalizing their simil-proofs, the process of formalization is not dissimilar from what mathematicians spend their time doing. In Jeremy Avigad's words,

> [t]he main moral I would like to extract […] is that formalization is continuous with the usual mathematical procedures for making claims and arguments precise. When one comes across a mathematical theorem, one may wonder whether the proof is correct, but once one convinces oneself that the proof is correct, there is no further question as to whether it can be formalized. Being correct means that one can supply details to any level of precision, down to the axioms and rules of a formal foundation if necessary.[14]

For our purposes, it is crucial to stress that this should hold for *all* areas of mathematics, even the ones in which we (or some among us) can simply *see* the results:[15]

---

[13] J. Avigad, "Reliability of mathematical inference," *Synthese,* 198: 7377-7399, 2021.

[14] J. Avigad, "A Formalization of the Mutilated Chessboard Problem," 2019, p. 5, available on-line at: https://www.andrew.cmu.edu/user/avigad/Papers/mutilated.pdf

[15] For a discussion on how visual reasoning in mathematics can be rigorous, see: S. De Toffoli, "Reconciling Rigor and Intuition," *Erkenntnis,* 86: 1783–1802, 2021.

> Spelling out spatial or visual intuitions in mathematical terms can be inordinately difficult, but we know how to do that, too, and the fact that we can do it is part and parcel of what we take such arguments to be properly mathematical. (*ibid.*)

This is precisely the norm that Enriques violates. Again, we cite Conforto:[16]

> Having once declared to him that I did not see the truth of a statement that he believed to be evident, but **which we had tried in vain to prove logically**, he stopped short […]

It is safe to assume that Enriques did not doubt the existence of a logical proof of the results he reached through his mathematical visions. Still, he was at times unable to find one – and this did not seem to bother him. But it might bother others. And we find ourselves sympathizing with his student Conforto and feeling the frustration he might have felt.

This is also a vivid illustration of the power conferred to mathematical arguments by what we called *rigor*. Rigorous arguments can be shared among rational, appropriately trained agents. The impotence of the Italian School is caused precisely by the communication impasses it faces due to the lack of rigor. But there is a caveat. If rigor is characterized functionally as what allows us to achieve reliability, then Enriques could argue that his practice was more rigorous than our putatively rigorous practice of creating detailed proofs that admit formalization. Let us return to the little dog. If we can access mathematical reality directly, *formal* rigor is unnecessary and may even hinder what really matters: getting things right. This might be an elitist standpoint, but this fact by itself does not delegitimize it as a viable position.

The goal of the present article is to investigate the details of the Italian School's *modus operandi* as well as its methodological reflections. We articulate different conceptions of mathematical rigor and link them to various conceptions of mathematical objectivity. The article proceeds as follows. In Section 2, we focus on how two of the leading members of the Italian School of algebraic geometry, Enriques and Severi, conceived of mathematical rigor. As it turns out, they did not understand rigor to be a unitary phenomenon. According to them, there are two distinct types of rigor. We show that from their perspective, which is almost paradoxical by contemporary standards, one type of rigor cannot be achieved without appeal to intuition. In Section 3, we outline how mathematical rigor can be linked to mathematical objectivity. Following Enriques and Severi's mathematical practice and their reflections thereon, we distinguish between two notions of objectivity: (1) objectivity as faithfulness to the facts and (2) objectivity as intersubjectivity.

---

[16] F. Conforto, "Federigo Enriques," art. cit.

## §2. SMALL-SCALE RIGOR AND LARGE-SCALE RIGOR

Throughout his entire mathematical life, Castelnuovo prized rigor. In particular, he was careful to distinguish between proofs (or better, simil-proofs) from mere plausibility arguments. Still, he employed the latter as cornerstones in some of his mathematical constructions. For instance, this is the case for his numerous applications of a number conservation principle:

> We must recognize that in establishing this concept **we rely more on intuition (and on various verifications) than on real mathematical reasoning**. Perhaps we can arrive at the proof by considering the curve in a higher-dimensional space as a partial intersection of several manifolds and by treating the problem of secant spaces algebraically; one would find that the number of solutions is independent of the particular position of the manifolds. But this kind of reasoning can only be done when the theory of curves in the higher-dimensional spaces is more complete. However, **we allow ourselves to take advantage of a not-yet-proven principle in order to solve a difficult problem, because we believe that, even** **with such attempts, it is to the benefit of science when we explicitly state what we admit and what we prove.**[17]

Several years later, in the comments added to the reprint of his memoir, issued on the occasion of his scientific jubilee, Castelnuovo remarks:

> This principle of degeneration [...] is simply admitted; the first proof that the splitting does not change the required numbers was given through topology (by exploiting Riemann surfaces) by F. Klein in 1892 [...]. From an algebraic perspective, one needs to show that the reducible curve may be obtained as a limit of an irreducible curve varying in a continuous system, that is, under very mild assumptions, what F. Severi proved in [...] 1921.[18]

While, for Castelnuovo, mathematical rigor seems to come after a more intuitive development of a problem, for Enriques rigor cannot be sharply separated from intuition:

> Even the usual question of whether mathematics should teach intuition or logic is vitiated by an imperfect vision of the value of teaching. As a matter of fact, **the presupposition of this question is that logic and intuition allow themselves to be separated as distinct intellectual faculties, whereas they are in fact two inseparable aspects of the same active process, which refer to each other.**[19]

---

In this context, Enriques uses the term *logic* as we use the term *rigor*. He distinguishes between two different notions falling under the same label:

> I need to state that logic encompasses more aspects than are usually seen by mathematics teachers. There is […] a **small-scale logic and a large-scale logic**: I mean the refined analysis of the exact thought process (almost the microscopic view of the elements that form the fabric of science), and, on the other hand, the study of the organic connections of the system – that is, the macroscopic view of science. **Now I fear that, for our mathematical teachers, small-scale logic holds too great a place in comparison to large-scale logic!** […] It is useless to develop the series of theorems of Euclidean geometry with impeccable deduction if we do not return to contemplating the built edifice, inviting the pupils to discern the truly significant geometric properties (e.g., the sum of the angles of a triangle and the Pythagorean theorem) from those having value only as links in the chain. That kind of democratic equality which some masters claim to establish among the demonstrated propositions – on the pretext that *everything is important* and therefore oblivion for nothing can be tolerated – only succeeds in deforming intelligences by depriving them of the light of evaluation, so that – borrowing words of a philosopher – the science that offers itself to the scholar in this way could be called "the infinite night, in which all the cows are black."[20]

We can characterize *small-scale rigor* as what is nowadays simply called *rigor*; it is concerned with the possibility of decomposing mathematical inferences into smaller inferences, and eventually reaching some basic inferences. *Large-scale rigor* is, however, a very different notion. It pertains to the large-scale structure of arguments and theories.

---

[20] *Ibid.*, pp. 9-10.

A similar distinction is also drawn by Severi. He ascribes to classical algebraic geometry a *substantial rigor*, meaning an essential faithfulness to the geometrical datum. This type of rigor is opposed to *formal rigor*, which refers to a strict axiomatic approach. In Severi's own words:

> The synthetic-intuitive turn of Poncelet's work was struck by the features of Cauchy, who judged the principle of continuity only as "a strong induction." The principle is, in fact, only an attempt at a general formulation and broad application of the theory of continuity of algebraic functions. If it had been possible at that time to give a precise statement of the principle, as an expression of an intuitive fact to be proved later, the use would have been quite legitimate because the principle would have become an axiom with an effective field applicability, as demonstrated by the very just applications Poncelet made of it. **We would thus have respected the rigor, which, according to Peano's witty phrase, is completely satisfied when we say, as in the oath to the tribunals, the whole truth and nothing but the truth: a rigor that I could call substantial to distinguish it from the formal rigor that each time calls for an axiomatic systematization from the letter A to the letter Ω.**[21]

On this occasion, Severi is addressing a French audience, which would be well aware of Cauchy's sharp critiques of Poncelet's incomplete proofs of correct results. Substantial rigor can be assimilated to Enriques's large-scale rigor and formal rigor to Enriques's small-scale rigor.

We have been exploring the possibility that rigor is what underwrites the reliability of mathematical arguments. A brief remark is in order. In the context of the present essay, we prefer talking about *(mathematical) arguments* rather than *proofs* or *simil-proofs*. This is because (i) arguments, like simil-proofs (but unlike proofs), might be incorrect, and (ii) arguments, unlike simil-proofs (and unlike proofs), might not be put forward as rigorous proofs, but simply as speculative reasons – such arguments were, for example, explicitly appealed to by Castelnuovo, as the quote above makes clear. For a mathematical argument to be reliable means that it is likely to be correct. But how do we assess such likelihood? By grasping how such an argument supports (if at all) its conclusion. It follows that the rigor of a mathematical argument can be conceptualized as what underwrites the possibility of such a grasping. This way of analyzing rigor is particularly useful in this context because grasping can also be conceived as a two-fold notion; consequently, the two types of rigor can be associated with two types of grasping:

1. **a local, step by step, grasping** of how (if at all) the [argument] supports its conclusion, plus the acceptance of its premises. Each step is recognized as valid, that is, each inferred item of information is held to be a logical consequence of earlier items of the simil-proof.[22]

2. **a holistic grasping** of how (if at all) the [argument] supports its conclusion, plus the acceptance of its premises. Such grasping often involves being aware of the structure of the argument at a large scale or how the argument and its major lemmas relate to what is already known about relevant mathematical areas, and so on.[23]

---

[21] F. Severi, "La Géometrie Algébrique Italienne. Sa Riguer, Ses Méthodes, Sès Problèmes," in *Opere Matematiche. Memorie e Note*, VI:9-55 [1-39], 1949, p. 2.
Our translation from French.

If rigor guarantees the reliability of arguments, it is plausible to think that its function is to allow the relevant audience to grasp how an argument supports (if at all) its conclusion. Small-scale rigor underwrites the possibility of local grasping, while large-scale rigor underwrites the possibility of holistic grasping.

Although these two types of grasping generally coexist, they are distinct. Differently than mathematical arguments in general, simil-proofs are put forward as rigorous proofs and therefore should make local grasping possible. It is a fact, however, that our practice of checking and evaluating simil-proofs does not always include a complete local grasping.[24] Different simil-proofs embody different types of rigor and therefore enable different forms of grasping.

---

[22] For our present purposes, we do not need to spell out what accepting the premises of a mathematical argument amounts to.
[23] S. De Toffoli, *Epistemic Roles of Mathematical Diagrams*, PhD Thesis, Stanford University, 2018, pp. 24-25.

Formal simil-proofs or highly detailed informal simil-proofs, for example, might make local grasping routine but also make holistic grasping quite hard or even practically impossible. This might seem at odds with the idea that small-scale rigor is what underwrites the reliability of simil-proofs. John Burgess remarks:

> A purported proof is much more likely to contain difficult-to-detect irreparable errors and so does not really amount to a genuine proof, if it is presented simply as a sequence of deductive steps, without enough hints as to overall strategy. This is why, for an expert, what looks to a non-expert like a mere outline of a proof may be more compelling than a lengthier account with a lot of fussy details, if the fussy details obscure the overall strategy.[25]

Setting aside automated checking of simil-proofs and focusing exclusively on humans checking simil-proofs, we see that excessive details may be harmful because they tend to hinder our holistic grasping:

> a purported proof that includes no hint of the overall strategy of argument, nor clues as to how it was discovered, and contains perhaps instead an excessive amount of fussy deductive detail better left to the reader […] may be a rigorous proof, containing no gaps or fallacies, but without a map of the woods we will have difficulty convincing ourselves that it is.[26]

---

[24] D. Fallis, "Intentional Gaps in Mathematical Proofs," *Synthese,* 134 (45–69), 2003.
[25] J. Burgess, *Rigor and Structure*, Oxford, Oxford University Press, 2015, p. 95.
[26] *Ibid.* p. 96.

A simil-proof overburdened with too many details can reduce its reliability.[27] However, for obvious reasons, too few details can be harmful as well. Holistic grasping alone does not guarantee that every detail in a simil-proof will be correct. Moreover, if the details cannot be supplied, we risk a communication impasse. That is precisely what happened when Enriques was unable to convince Conforto of the truth of his result. The fact that arguments published by the members of the Italian School did not always allow for local grasping can also be used to explain why errors in the published literature remained often unnoticed and why the Italian School ended when it did.

We are exploring the idea that rigor is what allows the relevant audience to grasp how a simil-proof supports (if at all) its conclusion. This idea allowed us to connect different types of rigor to different types of grasping. Note that the view on offer ties rigor to reliability but challenges the orthodoxy according to which rigor implies perfect reliability, or correctness. In fact, under the present assumptions even simil-proof that turn out to be fallacious can be rigorous (in both senses).

## §3. VARIETIES OF MATHEMATICAL OBJECTIVITY

In this section, we turn to the notion of mathematical objectivity. We submit that mathematical objectivity can be understood in relation to rigor and that, consequently, it can be understood as a two-fold notion. But it will be first useful to consider how scientific objectivity can be characterized in general.

Julian Reiss and Jan Sprenger[28] distinguish between three different characterizations of scientific objectivity:

---

[27] Here we are glossing over the intricate issue of how to identify simil-proofs. One might think that adding or subtracting details would not affect the simil-proof but only a presentation thereof. In the present context, we do not need to settle this difficult matter.
[28] J. Reiss and J. Sprenger, "Scientific Objectivity," *Stanford Encyclopedia of Philosophy*, no. Winter 2020: 1–72, 2020.

(1) faithfulness to the facts,

(2) absence of personal bias, and

(3) absence of normative commitments and value-freedom.[29]

In the context at hand, we are interested in the first two characterizations. This is because we aim to investigate the objectivity of mathematical arguments. Reiss and Sprenger's third characterization of objectivity has to do with changes in theories within science and is therefore irrelevant in the present context.

Let us now turn to the analysis of the objectivity of mathematical arguments interpreted either as some property that leads us faithfully to mathematical facts, or as the property of being neither subjective nor personal – that is, of being intersubjective.[30]

## 3.1 OBJECTIVITY AS FAITHFULNESS TO THE FACTS

The first form of objectivity assumes that science aims at correctly describing facts:

> The philosophical rationale underlying this conception of objectivity is the view that there are facts "out there" in the world and that it is the task of scientists to discover, analyze, and systematize these facts. "Objective" then becomes a success word: if a claim is objective, it correctly describes some aspect of the world.[31]

One may doubt whether it is possible to apply such a conception of objectivity to mathematics. After all, it is hardly straightforward to articulate what it means for mathematical facts to be "out there."[32] However, from the point of view of the Italian School of algebraic geometry, this conception of objectivity suits mathematics well. As we have seen, Castelnuovo, Enriques, and other members of their school conceived of algebraic surfaces in a quasi-concrete way, and their goal was to describe them faithfully using what they thought were reliable quasi-perceptual faculties.

---

[29] For expository purposes, we inverted the order of the last two items of the list.
[30] For a discussion of intersubjectivity in the context of an analysis of the nature of epistemic justification in general, see S. De Toffoli "Intersubjective Propositional Justification," in P. Silva Jr. and L. R G Oliveira (eds.), *Propositional and Doxastic Justification: New Essays on Their Nature and Significance*, Abingdon, Routledge, 2022.
[31] J. Reiss and J. Sprenger, "Scientific Objectivity," art cit., p.4.

From Enriques's perspective, it is large-scale rigor which allows mathematicians to better describe mathematical facts faithfully. This is because large-scale rigor enables a form of holistic grasping. An impressionistic argument might be better suited to guide (some of) us to see the mathematical facts in question over a detailed argument specifying each logical step. A couple of brushstrokes from a good painter can suffice to represent a little dog faithfully – we might not need the pixel-by-pixel description if we just want to be sure of the dog's existence.

All this is fine and good. But not for everyone. As we have seen, mathematical visualizations are not easily communicated. Consider the following piece of mathematical folklore, narrated by Fields medalist Vaughan Jones:

> Once, at a seminar, one of the world's best low-dimensional topologists was presenting a result. At a certain point another distinguished topologist in the audience intervened to say he did not understand how the speaker did a certain thing. The speaker gave an anguished look and gazed at the ceiling for at least a minute. The member of the audience then

---

[32] According to some, however, the exact opposite is true:

> Mathematicians may have even better experimental access to mathematical reality than the laboratory sciences have to physical reality. This is the point of modeling: a physical phenomenon is approximated by a mathematical model; then the model is studied precisely because it is more accessible. This accessibility also has had consequences for mathematics on a social level. Mathematics is much more finely subdivided into subdisciplines than physics, because the methods have permitted a deeper penetration into the subject matter.

A. Jaffe, and F. Quinn, "''Theoretical Mathematics'': Toward a Cultural Synthesis of Mathematics and Theoretical Physics," *Bulletin of the American Mathematical Society,* 29 (1): 1–13, 1993, p. 2.

affirmed "Oh, yes, I hadn't thought of that!" Visibly relieved, the speaker went on with his talk, glad to have communicated this point to the audience.[33]

Like Italian style algebraic geometry, low-dimensional topology is an area of mathematics in which heavy reliance on visualizations is common.

Again, the problem is that it should be possible to re-cast these visualizations as detailed arguments involving well-defined concepts and operations. It is only in this way that small-scale rigor and intuition can be reconciled. While this is often achieved in modern mathematical practice, there are notable exceptions. The work of the Italian School is one such exception – it is due to its failure to connect large-scale rigor to small-scale rigor that it was a problematic mathematical community. For this reason, it could only last for a relatively short period of time before succumbing to the new methodological guidelines of the French algebraic geometers from the 1940s to the 1960s, after which it only began being appreciated again in the mid-1960s. This is largely due to the efforts of David Mumford, who inaugurated the era of *neoclassical algebraic geometry* by applying Grothendieck's techniques to Enriques's legacy.[34]

A mathematical practice can only flourish if the methods it adopts are *shareable,* that is, if it can be communicated and used by different mathematicians. A mathematical practice which can be carried out only by exceptionally talented subjects, however, cannot live long. It is on this point that a separate conception of objectivity proves quite useful.

## 3.2 OBJECTIVITY AS INTERSUBJECTIVITY

Scientific objectivity can also be understood as a form of intersubjectivity —as freedom from personal biases. According to this view, science is objective to the extent that personal biases are absent from scientific reasoning, or that they can be eliminated in a social process. Perhaps all science is necessarily perspectival. Perhaps we cannot sensibly draw scientific inferences without a host of background assumptions, which may include assumptions about values. Perhaps all scientists are biased in some way. But objective scientific results do not, or so the argument goes, depend on researchers' personal preferences or experiences—they are the result of a process where individual biases are gradually filtered out and replaced by **agreed upon evidence**.[35]

---

[33] V. F. R. Jones, "A Credo of Sorts," in H.G. Dales and G. Oliveri (eds.), *Truth in Mathematics*, 203–14. Oxford University Press, 1998, p. 214.
[34] D. Mumford, *Lectures on Curves on an Algebraic Surface*, Princeton University Press, 1966 and D. Mumford, "Intuition and Rigor and Enriques's Quest," *Notices of the American Mathematical Society,* 58: 250–60, 2011.

We must qualify this characterization to make it applicable to the case at hand. Mathematical arguments constitute *agreed upon evidence* only when they are *shareable*:

An argument is *shareable* if its content and supposed correctness could be grasped by relevantly trained human minds from a (possibly enthymematic) perceptible instance of a presentation of it.[36]

As Conforto's recollections make it clear, Enriques's visions cannot constitute shareable arguments. That is why they cannot be "agreed-upon evidence." However, the failure of communication is not due to specific biases or idiosyncratic values. Instead, it is due to the appeal to faculties that are not part of our shared cognitive makeup.

In the case of science, measurements supply evidence all subjects can agree upon. In 1883 Lord Kelvin wrote:

---

[35] J. Reiss and J. Sprenger, "Scientific Objectivity," art cit., p. 28.
[36] S. De Toffoli, "Groundwork for a Fallibilist Account of Mathematics," art.cit., p. 830.

When you cannot express it in numbers, your knowledge is of a meagre and unsatisfactory kind; it may be the beginning of knowledge, but you have scarcely, in your thoughts, advanced to the stage of science, whatever the matter may be.[37]

In mathematics, evidence all subjects can agree upon is not always about numbers but is also about the logical structure of arguments. Objectivity as intersubjectivity is therefore achieved through small-scale rigor. A simil-proof has small-scale rigor if and only if it enables local, step-by-step grasping.[38] This type of grasping, at its extreme, does not require a subject at all – at that point, it becomes mechanical.

Large-scale rigor is usually more satisfactory than small-scale rigor. It is, for instance, associated with a phenomenology of "enlightenment."[39] However, it is not something that can always be achieved. More importantly, even when it *can* be achieved, it can't be achieved by everyone. Mathematical arguments lacking small-scale rigor risk a failure to be objective in this sense; that is, they risk a failure of being intersubjective.[40]

## §4. CODA

---

[37] Quoted in J. Reiss and J. Sprenger, "Scientific Objectivity," art cit., p. 29.

[38] Note that when we talk about small-scale rigor, we talk about simil-proofs and not about arguments in general. This is because aspiring to small-scale rigor is a characteristic of proofs and simil-proofs not shared with mathematical arguments in general.

[39] This phenomenon is described in: G.C. Rota, "The phenomenology of mathematical proof," in *Indiscrete Thoughts*, New York, Modern Birkhäuser Classics**,** pp. 134-150, 1997.

[40] As pointed out by an anonymous referee, in Enriques's wide philosophical production there are also passages suggesting that geometric intuition can be conceived as intersubjective. In particular, the psychological interpretation of logic defended in *Problems of Science* (F. Enriques, *Problemi della Scienza*, Bologna, Zanichelli, 1906) seems to provide a basis for the idea that geometric intuition is something that should be trained and is potentially available to everyone. For a similar suggestion, see also H. Poincaré, "La logique et l'intuition," *L'enseignement mathématique, 1*(5): 157-162, 1889.

Castelnuovo, Enriques, and Severi considered their mathematical practice to be a quasi-perceptual mental exploration of abstract mathematical objects existing independently of human subjects.[41] Their milestone result, the classification of algebraic surfaces, was reached thanks to heavy reliance on mathematical *seeing*, a form of geometric intuition. Their methods fell short (by a long way) of modern standards of rigor. In fact, this is acknowledged frankly by Castelnuovo himself. In the introduction of his edition of Enriques's posthumous monograph on algebraic surfaces, he writes:[42]

> The author himself takes care to warn right from the preface that the treatise, rather than expounding an already static and crystallized doctrine, aspires to arouse in the reader the desire to bring additions and improvements to various theories. And where the ground is less solid, the author warns the scholar.[43]

Nonetheless, they managed to produce important and lasting results because they satisfied a different type of rigor. It is however important to remember that this methodology, especially when not explicitly acknowledged, also led to the publications of incorrect and gappy proofs and even of false statements – this was particularly problematic in the last period of Severi's career.[44]

It is certainly true that the Italian School produced arguments that lacked small-scale rigor, which is today a sine-qua-non condition for something to be an acceptable piece of mathematics. Those arguments, however, enjoyed a different type of rigor, which the members of the School deemed to be superior. Their large-scale rigor is what made them track better, more direct paths to mathematical facts. If objectivity is understood as faithfulness to the facts – then these arguments were objective. We cannot resist citing, again, Conforto's words:

> Having once declared to him that I did not see the truth of a statement, which he believed to be evident, but which we had tried in vain to prove logically, he stopped short […] he spun his stick, pointing it toward a little dog on a window sill, saying to me: "Can't you see? For me it's as if you were telling me you don't see that little dog!"[45]

---

[41] To be sure, although Castelnuovo, Enriques and Severi shared a general outlook, their conception of mathematics was not uniform; see A. Brigaglia and C. Ciliberto, "Geometria Algebrica," in S. di Sieno, A. Guerraggio and P. Nastasi (eds.), *La Matematica Italiana Dopo l'Unità. Gli Anni Tra Le Due Guerre Mondiali*, 185–320, Marcos y Marcos, 1998.
[42] Enriques died suddenly in 1946.
[43] F. Enriques, *Le Superfici Algebriche, op. cit.,* p. vi.
[44] See, for example, C. Ciliberto, C. and E. Sallent Del Colombo, "Francesco Severi: il suo pensiero matematico e politico prima e dopo la Grande Guerra," Preprint, ArXiv: 1807.05769, 2018.

Conforto did not see – but this was *his* problem. For the visionary Enriques, small-scale rigor would almost be a hindrance – a way of decomposing reality that, like a printer, annihilates all meaning. However, without it there is no way for the mathematical community to systematically check the correctness of the simil-proofs it produces – and thus there is a significant danger of unnoticed mistakes in the published literature.

Faithfulness to the facts may not be a something everyone can achieve, but it can be reached by a "powerfully intuitive spirit." It is because of this that such a notion of objectivity is in tension with objectivity *as* intersubjectivity. While the latter is democratizing, the former remains elitist. Enriques himself, once said: "intuition is the aristocratic way of discovery, rigour the plebian way."[46]

Objectivity as faithfulness to mathematical facts can lead to substantial advances but is inadequate and dangerous in the end. As we see it, the objectivity for which mathematicians strive today integrates both of these conceptions: it is a matter of faithfully describing mathematical reality *and* of doing so in a way that can be appreciated by any rational subject with the appropriate background training.

---

[45] F. Conforto, "Federigo Enriques," art. cit.
[46] W. V. D. Hodge. "Federigo Enriques", *Hon. F. R. S. E., Royal Society of Edinburgh Year Book, 1948/49*, 13–16. Quoted in D. Babbitt, and J. Goodstein, "Federigo Enriques's Quest to Prove the 'Completeness Theorem'," *Notices of the American Mathematical Society*, 58: 240-249, 2011, p. 240.